\documentclass{tac}
\usepackage{xypic}
\usepackage{enumerate,amsmath,amssymb,latexsym,QED}
\usepackage[ps2pdf, bookmarks=true,
            bookmarksnumbered=true, breaklinks=true,
            pdfstartview=FitH, hyperfigures=false,
            plainpages=false, naturalnames=true,
            colorlinks=true,
            pdfpagelabels]{hyperref}
\usepackage[usenames]{color}

\xyoption{curve}

\newtheoremrm{example}{Example}
\newtheoremrm{defn}{Definition}

\newtheoremrm{Def}{Definition}

\newtheorem{prop}{Proposition}
\newtheorem{thm}[example]{Theorem}
\newtheoremrm{rem}{Remark}
\newtheorem{cor}{Corollary}
\newtheorem{lem}{Lemma}

\newcommand{\colim}{\operatorname{colim}}

\def\<{\langle}
\def\>{\rangle}
\def\Z{\mathbb{Z}}

\def\X{\mathsf{X}}
\def\B{\mathsf{B}}
\def\xybiglabels{\def\labelstyle{\textstyle}}
\def\id{\operatorname{id}}

\def\Gpds{\mathsf{Gpd}}

\def\Set{\mathsf{Set}}

\newcommand{\Sets}{\mathsf{Set}}

\newcommand{\XMod}{\mathsf{XMod}}

\def\Crs{\mathsf{Crs}}

\def\M{\mathsf{M}}

\newcommand{\eps}{\varepsilon}

\newcommand{\io}{^{-1}}
\newcommand{\C}{\mathsf{C}}

\def\geq{\geqslant}

\def\leq{\leqslant}

\def\QED{\hfill $\Box$}

\newcommand{\labto}[1]{\stackrel{#1}{\longrightarrow}}

\def\phi{\varphi}

\newcommand{\I}{{\mathcal I}}

\def\Gpd{\mathsf{Gpd}}

\newcommand{\Ob}{\operatorname{Ob}}

\def\Mod{\mathsf{Mod}}
\def\XSq{\mathsf{XSq}}

\def\epsilon{\varepsilon}
\def\P{\Phi}
\def\D{\mathsf{D}}

\def\R{\mathbb{R}}
\def\Ob{\operatorname{Ob}}
\def\a{\alpha}
\def\b{\beta}
\def\2-Gpd{\mathsf{2\mbox{-}Gpd}}
\begin{document}

\author{Ronald Brown and Rafael Sivera}
\thanks{This work was partially supported by a Leverhulme Emeritus Fellowship 2002-4 for
the first author.}


\address{School of Computer Science \\
Bangor University,  \\ Gwynedd LL57 1UT, \\  U.K.\\[5pt]
Departamento de {Geometr\'{i}a y Topolog\'{i}a}\\
          Universitat de Val\`encia,\\ 46100 Burjassot, Valencia\\
         Spain}
\title{Algebraic colimit calculations in homotopy theory  \\ using fibred and cofibred categories}
\eaddress{r.brown@bangor.ac.uk, \CR Rafael.Sivera@uv.es}

         \maketitle
\begin{center}
  Bangor University Maths Preprint 06.08
\end{center}

\begin{abstract}
{Higher Homotopy van Kampen Theorems allow the computation as
colimits of certain homotopical invariants of glued spaces. One
corollary is to describe homotopical excision in critical dimensions
in terms of induced modules and crossed modules over groupoids. This
paper shows how fibred and cofibred categories give an overall
context for discussing and computing such constructions, allowing
one result to cover many cases.  A useful general result is that the
inclusion of a fibre of a fibred category preserves connected
colimits. The main homotopical application are to pairs of spaces
with several base points, but we also describe briefly the situation
for triads. }\footnote{29/09/08}
\end{abstract}

\tableofcontents \addtocontents{References}

\section{Introduction} \label{sec:intro}
One of our aims is to give the  framework of fibred and cofibred
categories to certain  colimit calculations of algebraic homotopical
invariants for spaces with parts glued together: the data here is
information on the invariants of the parts, and of the gluing
process.

The second aim is to advertise the possibility of such calculations,
which are based on various Higher Homotopy van Kampen Theorems
(HHvKTs)\footnote{This name for certain generalisations of van
Kampen's Theorem was recently suggested by Jim Stasheff. }, proved
in 1978-87. These are of the form that a {\it homotopically defined}
functor
$$\Pi: (\text{topological  data}) \to (\text{algebraic  data}) $$
preserves certain colimits \cite{BH78:sec,BH81:col,BL87}, and where
the algebraic data contains non-Abelian information.

The antecedent for dimension 1 of such functors as $\Pi$  is the
fundamental group $\pi_1(X,a)$ of based spaces:   the van Kampen
theorem in its form due to Crowell, \cite{Cr59}, gives that this
$\pi_1$ preserves certain pushouts.

The next advance was the fundamental groupoid functor $\pi_1(X,A)$
from spaces $X$ with a {\it set} $A$ of base points to groupoids:
groupoids in effect carry information in dimensions 0 and 1. For
example, the groupoid van Kampen theorem \cite{B67,B2006} gives the
fundamental group of the circle as the infinite cyclic group $\C$
obtained in the category of groupoids from the (finite) groupoid $\I
\cong \pi_1([0,1],\{0,1\})$ by identifying $0,1$. That is, we have
the analogous pushouts
\begin{equation}\label{eq:circle1}
 \vcenter{\xymatrix{\{0,1\} \ar[r] \ar[d] &\ar [d]  \{0\}\\
  [0,1] \ar [r]& S^1 }} \overset{\pi_1}{\mapsto}
   \vcenter{\xymatrix{\{0,1\} \ar[r] \ar[d] &\ar [d]  \{0\}\\
  \I \ar [r]& \C }}
\end{equation}
the first in the category of spaces, the second in the category of
groupoids, and $\pi_1$ takes the first with the appropriate sets of
base points, to the second. The aim is to have a similar argument at
the module level for $\pi_n, n>1$.

One of the problems of obtaining such results in higher homotopy
theory is that low dimensional identifications of spaces usually
affect high dimensional homotopy invariants. To cope with this fact,
the algebraic data for the values of the functor $\Pi$ must have
structure interacting from low to high  dimensions, in order to
model how the spaces are glued together.

An example of such dimensional interaction is the operation  of the
fundamental group on the higher homotopy groups,  which fascinated
the early workers in homotopy theory  (private communication, J.H.C.
Whitehead, 1957).  This operation can be seen as  necessitated by
the dependence of these groups on the base point, but  is somewhat
neglected in algebraic topology, perhaps for perceived  deficiencies
in modes of calculations.

A recognised reason for the difficulty of traditional invariants
such as homotopy groups in dealing with gluing spaces, and so
obtaining colimit calculations, is the failure of excision. If $X$
is the union of open sets $U,V$ with intersection $W$ then the
inclusion of pairs
$$\varepsilon:(V,W) \to (X,U)$$ does not induce an isomorphism
of relative homotopy groups as it does for singular homology.
However there is a residue of homotopical excision which  can be
deduced from a  HHvKT, as in \cite{BH81:col}; it requires
connectivity conditions and gives information in only the critical
dimension, but this information involves the additional structure of
operations. We state it at this stage for a single base point, as
follows:

\noindent  {\bf Homotopical Excision Theorem (HET)} {\it If $U,V$
and  $W=U \cap V$ are path connected with base point $a \in W$, and
$(V,W)$ is $(n-1)$-connected, $n \geq 3$, then $(X,U)$ is
$(n-1)$-connected and the excision map
$$\varepsilon_*:\pi_n(V,W,a) \to \pi_n(X,U,a)$$ presents the $\pi_1(U,a)$-module $\pi_n(X,U,a)$
as induced  by the morphism of fundamental groups $\pi_1(W,a)\to
\pi_1(U,a)$ from the  $\pi_1(W,a)$-module  $\pi_n(V,W,a)$. The same
holds for $n=2$ with `module' replaced by `crossed module'. }

We show that this theorem implies the classical Relative Hurewicz
Theorem, (Corollary \ref{cor:RHT}) for which in more usual  terms
see for example \cite{WG78}. The case $n=2$ implies a deep theorem
of Whitehead on free crossed modules (Corollary
\ref{cor:freeonattaching}) (of which \cite{B80} explains the
original proof), and indeed allows the more general calculation of
$\pi_2(X \cup CA,X,a)$ as a crossed $\pi_1(X,a)$-module (Corollary
\ref{cor:attachcone}). That result is applied in \cite{BW1,BW03} to
calculate crossed modules representing some homotopy 2-types of
mapping cones of maps of classifying spaces of groups.

The notion of {\it inducing} used here is both an indication of the
existence of some universal properties in higher homotopy theory,
and also of the tools needed to exploit such  properties as there
are. In the case $n>2$ the inducing construction is well known for
modules over groups: if $M$ is a module over the group $G$, and $f:
G \to H$ is a morphism of groups, then the induced module  $f_*(M)$
is isomorphic to $M \otimes _{\Z G} \Z H$. For $n=2$, the analogous
construction is non-Abelian.

We give the approach of fibred and cofibred categories to this
inducing construction.

Also the following simple example of algebraic modelling of a gluing
process suggests a further generalisation is needed. Consider the
following maps involving the $n$-sphere $S^n$, the first an
inclusion and the second an identification:
\begin{equation}\label{eq:nsphere}
  S^n \labto{i} S^n \vee [0,1] \labto{p} S^n \vee S^1.
\end{equation}
Here $n \geq 2$, the base point {$\mathbf{0}$} of $S^n$ is in the
second space  identified with $0$, and the map $p$ identifies $0,1$
to give the circle $S^1$.

This clearly  requires an algebraic theory dealing with many base
points, and so a generalisation from groups to groupoids. The geometry
of the base points and their interconnections can then play a role in the theory and applications.

This programme has proved successful for  the algebraic data of {\it
crossed complexes over groupoids}, or equivalent structures, and for
{\it cat$^n$-groups}, or, equivalently, crossed $n$-cubes of groups.
In both categories there is a HHvKT asserting  that a homotopically
defined functor $\Pi$ preserves certain colimits, and this leads to
new calculations of homotopical invariants.

The calculation of colimits in these algebraic categories usually
requires working up in dimensions.  So it is useful to consider for
these hierarchical structures the forgetful functors from higher
complex  structures to lower structures.

The prototype again is groupoids, with the object functor from
groupoids to sets. The fibre of this functor over a set $I$ is
written $\Gpd_I$: the objects of this category are groupoids with
object set $I$ and morphisms are functors which are the identity on
$I$. This fibre has nice properties, and is what is called {\it
protomodular};  in rough terms it has properties analogous to those
of the category of groups. In particular, there is a notion of
normal subgroupoid, and of coproduct, so that colimits can be
calculated in $\Gpd_I$ as quotients  of a coproduct by a normal
subobject.

However the whole category $\Gpd$ is of interest for homotopical
modelling. The functor $ \Ob: \Gpd \to \Set$ has the properties of
being a fibration and a cofibration in the sense of Grothendieck,
\cite{Gro68-cofibred}. So a function $u:I \to J$ of sets gives rise
to a pair of  functors $u^*: \Gpd_J \to \Gpd_I$, known often as
pullback, and $u_*: \Gpd_I \to \Gpd_J$, which could be called
`pushin' or `push forward', such that $u_*$ is left adjoint to
$u^*$. So for $G \in \Gpd_I$ we have a morphism $u': G \to u_*G$
over $u$ with a universal property which is traditionally  called
`cocartesian', and it is also said that $u_*G$ is the object {\it
induced from $X$ by $u$}.

It is interesting to see that the main parts of the pattern for
groupoids goes over to the general case. Our  categorical results
are for a fibration of categories $\Phi: \X \to \B$. The first and a
main result  (Theorem \ref{thm:preservecolimfibres}) is that the
inclusion $\X_I \to \X$ of a fibre $\X_I$ of $\Phi$ preserves
colimits of connected diagrams. Our second set of results relates
pushouts in $\X$  to the construction of the functors $u_*: \X_I \to
\X_J$ for $u$ a morphism in $\B$ (Proposition
\ref{prop:induced-pushout}). Finally, we show how the combination of
these results uses the computation of colimits in $\B$ and in each
$\X_I$ to give the computation of colimits in $\X$ (Theorem
\ref{thm:compcolim-fibred}).

These general results are developed in the spirit of `categories for
the working mathematician' in sections \ref{sec:fibcat} to
\ref{sec:pushouts}. We illustrate the use of these results for
homotopical calculations not only in groupoids (section
\ref{sec:groupoidsinduced}), but also for crossed modules (section
\ref{sec:xmodGpd}) and for modules, in both cases over groupoids.
Finally we give a brief account of some relevance to crossed squares
(section \ref{sec:crossedsquares}), as an indication of a more
extensive theory, and in which these ideas need development.

We are  grateful to Thomas Streicher for his Lecture Notes
\cite{Streicher-fibred},  on which the following account of fibred
categories is based, and for helpful comments leading to
improvements in the proofs. For further accounts of fibred and
cofibred categories, see
\cite{gray-fibred,Borceux-handb2,vistoli-fibredcat} and the
references there. The first paper gives analogies between fibrations
of categories and Hurewicz fibrations of spaces.

\section{Fibrations of categories}\label{sec:fibcat}
We recall the definition of fibration of categories.

\begin{Def}
Let $\Phi : \X \to \B$ be a functor. A morphism $\phi: Y \to X$ in
$\X$ over $u:=  \Phi (\phi)$ is called {\em cartesian}
\index{cartesian!morphism}  if and only if for all $v: K \to J$ in
$\B$ and $\theta: Z \to X$ with $\Phi (\theta)=u  v$ there is a
unique morphism $\psi: Z \to Y$ with $\Phi (\psi)=v$ and $\theta =
\phi  \psi$.

This is illustrated by the following diagram:

$$\xybiglabels\xymatrix@R=2.8pc@C=2.5pc{Z \ar @{..>}[r]_\psi \ar@/^1.5pc/ [rr] ^\theta
& Y \ar [r] _ \phi & X & \ar [d] ^ \Phi\\
K \ar [r]_v \ar @/^1.5pc/ [rr] ^{u v} & J \ar [r] _ u & I & }
$$
\hfill $\Box$
\end{Def}

It is straightforward to check that cartesian morphisms are closed
under composition, and that $\phi$ is an isomorphism if and only if
$\phi$ is a cartesian morphism over an isomorphism.

A morphism $\alpha: Z \to Y$ is called {\em vertical} (with
respect to $\Phi$) \index{vertical} if and only if $\Phi (\alpha)$
is an identity morphism in $\B$. In particular, for $I \in \B$ we
write $\X_I$,  called the {\it fibre over $I$}, \index{fibre!of a
functor} for the subcategory of $\X$ consisting of those morphisms
$\alpha$ with $\Phi (\alpha)= \id_I$.
\begin{Def} The functor
$\Phi : \X \to \B$ is a {\em fibration} \index{fibration!functor} or
{\em category fibred over $\B$} \index{category!fibred} if and only
if for all $u: J \to I$ in $\B$ and $X \in \X_I$ there is a
cartesian morphism $\phi: Y \to X$ over $u$: such a $\phi$ is called
a {\it cartesian lifting} \index{cartesian lifting} of $X$ along
$u$.
 \rule{0em}{1ex} \hfill $\Box$
\end{Def}

Notice that cartesian liftings of $X \in \X_I$ along $u:  J \to I$
are unique up to vertical isomorphism: if $\phi: Y \to X$ and
$\psi: Z \to X$ are cartesian over $u$, then there exist vertical
arrows $\alpha: Z \to Y$ and $\beta: Y \to Z$ with $\phi  \alpha=
\psi$ and $\psi  \beta= \phi$ respectively, from which it follows
by cartesianness of $\phi$ and $\psi$ that $\alpha \beta= id_Y$
and $\beta \alpha=\id _Z$ as $\psi \beta\alpha=\phi \alpha= \psi =
\psi \id_Y$ and similarly $\phi \beta\alpha= \phi \id_Y$.

\begin{example}\label{example:Ob}
The forgetful functor,  $\Ob:\Gpd \to \Sets$, from the category of
groupoids to the category of sets is a fibration. We can for a
groupoid $G$ over $I$ and function $u: J \to I$  define the
cartesian lifting $\phi: H \to G$ as follows: for $j,j' \in J$ set
$$H(j,j')= \{(j,g,j')\mid g \in G(uj,uj')\}$$ with composition
$$(j_1,g_1,j'_1)(j,g,j')= (j_1,g_1g,j'),$$  with $\phi$ given by
$\phi(j,g,j') = g$. The universal property is easily verified. The
groupoid $H$ is usually called the {\it pullback} of $G$ by $u$.
This is a well known construction (see for example \cite[\S
2.3]{Mackenzie2005}, where pullback by $u$ is written
$u^{\downarrow\downarrow})$. \rule{2em}{0ex} \QED
\end{example}

\begin{Def}\label{def:reindex}
If $\Phi :\X \to\B$ is a fibration, then using the axiom of choice
for classes we may select for every $u:J \to I$ in $\B$ and $X \in
\X_I$ a cartesian lifting of $X$ along $u$
$$u^X: u^* X \to X.$$
Such a choice of cartesian liftings is called a {\it cleavage} or
{\it splitting} \index{categories!fibration of!cleavage}
\index{categories!fibration of!splitting} of $\Phi $.

If we fix the morphism $u: J \to I$ in $\B$, the splitting gives a
so-called {\it reindexing functor} \index{categories!fibration
of!reindexing $u^*$}
$$u^*: \X_I \to \X_J$$ defined on objects by
$X \mapsto u^*X$ and the image of a morphism $\alpha : X \to Y$ is
$u^*\alpha$ the unique vertical arrow commuting the diagram:
$$\xybiglabels \xymatrix@C=6pc@R=3pc{u^*X \ar [r] ^{u^X} \ar@{-->}[d] _{u^*\alpha}
 & X \ar[d] ^\alpha \\
u^* Y \ar [r] _{u^Y}& Y}$$ \hfill $\Box$
\end{Def}
We can use this reindexing functor to get an adjoint situation for
each $u : J \to I$ in $\B$.

\begin{prop}\label{prop:fibadjoint} Suppose $\Phi: \X \to \B$ is a fibration of categories, $u
: J \to I$ in $\B$, and a reindexing functor $u^*: \X_I \to \X_J$ is
chosen. Then there is a bijection
$$ \X_J(Y, u^*X) \cong \X_u(Y,X)$$
natural in $Y \in \X_J, \, X \in \X_I$ where  $\X_u(Y,X)$ consists
of those morphisms $\alpha \in \X(Y,X)$ with $\Phi (\alpha)= u$.
\end{prop}
\begin{proof}
This is just a restatement of the universal properties concerned.
\end{proof}

In general for composable maps $u: J \to I$ and $v:K \to J$ in $\B$
it does not hold that $$v^*  u^* =(u  v)^* $$ {as may be seen with
the fibration of Example \ref{example:Ob}}. Nevertheless there is a
natural equivalence  $c_{u,v}:v^* u^* \simeq (u  v)^* $
as shown in the following diagram in which the full arrows are
cartesian and where $(c_{u,v})_X$ is the unique vertical arrow
making the diagram commute:
$$ \xybiglabels \xymatrix@C=7pc@R=3pc{ v^* u^* X \ar @{-->} [d] _{(c_{u,v})_X}^ \cong
\ar[r] ^{v^{u^*X}}& u^* X \ar[d] ^{u^X}\\
(u v)^* X \ar [r]  _{(u v)^X} & X }$$ Let us consider this
phenomenon for our main examples:

\begin{example}
1.- Typically, for $\Phi _\B= \partial_1: \B^2 \to \B$, the
fundamental fibration for a category with pullbacks, we do not know
how to choose pullbacks in a functorial way.

2.- In considering groupoids as a fibration over sets, if $u: J \to
I$ is a map, we have a reindexing functor, also called {\it
pullback}, $u^*: \Gpds_I \to \Gpds_J$. We notice that $v^* u^*Q$ is
naturally isomorphic to, but not identical to $(u  v)^* Q$. \QED
\end{example}

A result which aids understanding of  our calculation of pushouts
and some other colimits of groupoids, modules,  crossed complexes
and higher categories is the following. Recall that a category $\C$
is {\it connected } \index{connected!category}
\index{category!connected} if for any $c,c' \in C$ there is a
sequence of objects $c_0=c,c_1, \ldots,c_{n-1}, c_n=c'$ such that
for each $i=0, \ldots, n-1$ there is a morphism $c_i \to c_{i+1} $
or $c_{i+1}  \to c_i$ in $\C$. The sequence of morphisms arising in
this way   is called a {\it zig-zag} \index{zig-zag} from $c$ to
$c'$ of length $n$.

\begin{thm} \label{thm:preservecolimfibres}
Let $\Phi : \X \to \B $ be a fibration, and let $J \in \B$. Then the
inclusion $i_J: \X_J \to \X$ preserves colimits of connected
diagrams.
\end{thm}
\begin{proof}  We will need the following diagrams:
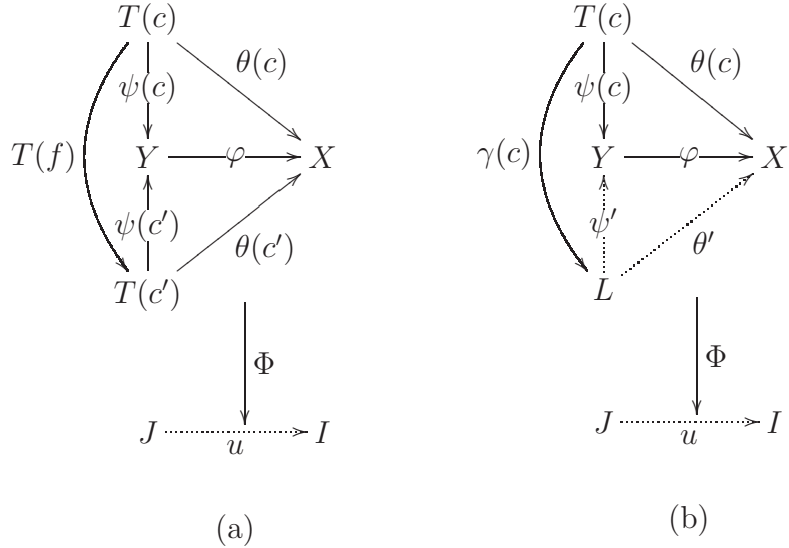
\begin{equation*}
\xybiglabels \xymatrix@R=3pc @C=1.5pc{ T(c) \ar@/_2pc/ [dd] _{T(f)}
\ar [d] |
{\psi(c)} \ar [drr] ^{\theta(c)} && \\
Y \ar [rr] | \phi && X \\
T(c')  \ar [u] | {\psi(c')} \ar [urr] _{\theta(c')} &\ar [d]^\Phi&\\
 J \ar @{..>} [rr]_u & &I\\
\ar @{}[rr]^{ \txt{(a)}} && }\qquad \qquad \xymatrix@R=3pc@C=1.5pc{
T(c) \ar@/_2pc/ [dd] _{\gamma(c)} \ar [d] |
{\psi(c)} \ar [drr] ^{\theta(c)} && \\
Y \ar [rr] | \phi && X \\
L  \ar @{..>}[u] | {\psi'} \ar @{..>}[urr] _{\theta'} &\ar [d]^\Phi&\\
 J \ar @{..>} [rr]_u & &I\\\ar @{}[rr]^{ \txt{(b)}} && }
\end{equation*}

Let $T: \C \to \X_J$ be a functor from a small connected category
$\C$ and suppose $T$ has a colimit $L \in \X_J $. So we regard $L$
as a constant functor $L:\C \to \X_J$ which  comes with a
universal cocone $\gamma: T \Rightarrow L$ in $\X_J$. Let $i_J:
\X_J \to \X$ be the inclusion. We prove that the natural
transformation $i_J \gamma:i_J T \Rightarrow i_J L$ is a colimit
cocone in $\X$.

We use the following lemma.

\begin{lem}
Let $X \in \X$, with $\Phi (X)=I$,  be  regarded as a constant
functor $X: \C \to \X$ and let $\theta: i_J T \Rightarrow  X$ be a
natural transformation, i.e. a cocone. Then \\ (i) for all $c \in
\C$, $u=\Phi (\theta(c)): J \to I$ in $\B$ is independent of $c$,
and\\ (ii) the  cartesian lifting $Y \to X$ of $u$ determines a
cocone $\psi: T \Rightarrow Y$.
\end{lem}
\begin{proof}
The natural transformation $\theta$  gives for each object the
morphism $\theta(c) : T(c) \to X$ in $\X$. Since $\C$ is connected,
induction on the length of a ziz-zag shows it is sufficient to prove
(i) when there is a morphism $f: c \to c'$ in $\C$. By naturality of
$\theta$, $\theta(c')T(f)= \theta(c)$. But $\Phi T(f)$ is the
identity, since $T$ has values in $\X_J$, and so $\Phi (\theta(c)) =
\Phi (\theta(c'))$. Write $u= \Phi (\theta(c))$.

Since $\Phi $ is a fibration, there is a  $Y \in \X_J$ and a
cartesian lifting $\phi: Y \to X$ of $u$.  Hence for each $c \in
\C$ there is a unique vertical morphism $\psi(c): T(c) \to Y$ in
$\X_J$ such that $\phi \psi(c) = \theta(c)$. We now prove that
$\psi$ is a natural transformation $T \Rightarrow Y$ in $\X_J$,
where $Y$ is regarded as a constant functor.

To this end, let $f: c \to c'$ be a morphism in $\X_J$. We need to
prove $\psi(c)= \psi(c') T(f)$.

The outer part of diagram (a) commutes, since $\theta$ is a natural
transformation. The upper and lower triangles commute, by
construction of $\phi$. Hence $$ \phi \psi(c)= \theta(c)=
\theta(c')T(f)= \phi \psi(c') T(f).$$ Now $T(f)$, $\psi(c)$ and $
\psi(c')$ are all vertical.  By the universal property of $\phi$,
$\psi(c)= \psi(c') T(f)$, i.e. the left hand cell commutes. That is,
$\psi$ is a natural transformation $T \Rightarrow Y$ in $\X_J$.
\end{proof}

To return to the theorem, since $L$ is a colimit in $\X_J$,  there
is a unique vertical morphism $\psi': L \to Y$ in the right hand
diagram (b) such that for all $c \in C$, $\psi'  \gamma(c) =
\psi(c)$. Let $\theta'= \phi  \psi': L \to X$. This gives a
morphism $\theta': L \rightarrow X$ such that $\theta' \gamma(c) =
\theta(c)$ for all $c$, and, again using universality of $\phi$,
this morphism is unique. \rule{1em}{0ex}
\end{proof}
\begin{rem}\label{rem:conn-colim}
The connectedness assumption is essential in the Theorem. Any small
category $\C$ is the disjoint union of its connected components. If
$T:\C \to \X$ is a functor, and $\X$ has colimits, then $\colim T$
is the coproduct (in $\X$) of the $\colim T_i$ where $T_i$ is the
restriction of $T$ to a component $\C_i$. But given two objects in
the same fibre of $\Phi: \X \to \B $, their coproduct in that fibre
is in general not the same as their coproduct in $\X$. For example,
the coproduct of two groups in the category of groups is the free
product of groups, while their coproduct as groupoids is their
disjoint union. \hfill $\Box$
\end{rem}
\begin{rem}
A common application of the theorem is that the inclusion $\X_J \to
\X$ {preserves} pushouts. This is relevant to our application of
pushouts in section \ref{sec:pushouts}. Pushouts are used to
construct free crossed modules as a special case of induced crossed
modules, \cite{BH78:sec}, and to construct free crossed complexes as
explained in \cite{BH91,BG89}. \QED
\end{rem}
\begin{rem}
George Janelidze has pointed out a short proof of Theorem
\ref{thm:preservecolimfibres} in the case $\Phi$ has a right
adjoint, and so preserves colimits, which applies to our main
examples here. If the image of $T$ is inside $\Phi(b)$, then $\Phi
T$ is the constant diagram whose value is $\{b, 1_b\}$, and if $\C$
is connected this implies that $\colim \Phi T = b$. But if
$\Phi(\colim T) = \colim \Phi T = b$, then $\colim T$ is inside
$\Phi(b)$. \QED
\end{rem}

\section{Cofibrations of categories}
We now give the duals of the above results.

\begin{Def}
Let $\Phi : \X \to \B$ be a functor. A morphism $\psi: Z \to Y$ in
$\X$ over $v:=  \Phi (\psi)$ is called {\em cocartesian}  if and
only if for all $u: J \to I$ in $\B$ and $\theta: Z \to X$ with
$\Phi (\theta)=u  v$ there is a unique morphism $\phi: Y \to X$
with $\Phi (\phi)=u$ and $\theta = \phi  \psi$.

This is illustrated by the following diagram:

$$\xybiglabels\xymatrix@R=2.8pc@C=2.5pc{Z \ar [r]_\psi \ar@/^1.5pc/ [rr] ^\theta
& Y \ar @{..>}[r] _ \phi & X & \ar [d] ^ \Phi\\
K \ar [r]_v \ar @/^1.5pc/ [rr] ^{u v} & J \ar [r] _ u & I & }
$$
\hfill $\Box$
\end{Def}

It is straightforward to check that cocartesian morphisms are
closed under composition, and that $\psi$ is an isomorphism if and
only if $\psi$ is a cocartesian morphism over an isomorphism.

\begin{Def} The functor
$\Phi : \X \to \B$ is a {\em cofibration}
\index{cofibration!functor} or {\em category cofibred over $\B$}
\index{category!cofibred} if and only if for all $v: K \to J$ in
$\B$ and $Z \in \X_K$ there is a cartesian morphism $\psi: Z \to Z'$
over $v$: such a $\psi$ is called a {\it cocartesian lifting}
\index{cocartesian lifting} of $Z$ along $v$.
 \rule{0em}{1ex} \hfill $\Box$
\end{Def}

The cocartesian liftings of $Z \in \X_K$ along $v:  K \to J$ are
also unique up to vertical isomorphism.

\begin{rem}
As in Definition \ref{def:reindex}, if $\Phi :\X \to\B$ is a
cofibration, then using the axiom of choice for classes we may
select for every $v:K \to J$ in $\B$ and $Z \in \X_K$ a
cocartesian lifting of $Z$ along $v$
$$v_Z: Z \to v_*Z.$$

{Under these conditions, the functor $v_*$ is commonly said to give
the objects {\it induced } \index{categories!cofibration of!functor
induced $v_*$} by $v$.} Examples of induced crossed modules of
groups are developed in \cite{BW03}, following on from the first
definition of these in \cite{BH78:sec}. \hfill $\Box$
\end{rem}

We now have the dual of Proposition \ref{prop:fibadjoint}.

\begin{prop} \label{prop:cofibadjoint}
For a cofibration $\Phi: \X \to \B$, a choice of cocartesian
liftings of $v : K \to J $ in $\B$ yields a functor $v_*: \X_K \to
\X_J$, 
and an adjointness
$$\X_J(v_*Z,Y) \cong \X_v(Z,Y)  $$ for all $Y \in \X_J,
\, Z \in \X_K$.
\end{prop}

We now state the dual of Theorem \ref{thm:preservecolimfibres}.

\begin{thm}
Let $\Phi: \X \to \B$ be a category cofibred over $\B$. Then the
inclusion of each fibre of $\Phi$ into $\X$ preserves limits of
connected diagrams.
\end{thm}

Many of the examples we are interested in are both fibred and
cofibred. For them we have an adjoint situation.

\begin{prop}
For a functor  $\Phi: \X \to \B$ which is both a fibration and
cofibration, and a morphism $u:J  \to I$ in $\B$, a choice of
cartesian and cocartesian liftings of $u$ gives an adjointness $$
\X_J(Y, u^*X) \cong \X_I(u_*Y,X)$$ for $Y \in \X_J, \, X \in
\X_I$.
\end{prop}

It is interesting to get a characterisation of the cofibration
property for a functor that already is a fibration. The following is
a useful weakening of the condition for cocartesian in the case of a
fibration of categories.

\begin{prop}\label{prop:cocartesianweakening}
Let $\Phi: \X \to \B $ be a fibration of categories. Then $\psi:Z
\to Y$ in $\X$ over $v:K \to J$ in $\B$ is cocartesian if only if
for all $\theta':Z \to X'$ over $v$ there is a unique morphism
$\psi': Y \to X'$ in $\X_J$ with $\theta'=\psi' \psi$.
\end{prop}
\begin{proof}
The `only if' part is trivial. So to prove `if' we have to prove
that for any $u : J \to I$ and $\theta : Z \to X$ such that
$\Phi(\theta)=uv$, there exists a unique $\phi : Y \to X$ over $u$
completing the diagram

$$\xybiglabels\xymatrix@R=2.8pc@C=2.5pc{Z \ar [r]_\psi
\ar@/^1.5pc/ [rr] ^\theta
& Y \ar @{..>}[r] _ \phi & X & \ar [d] ^ \Phi\\
K \ar [r]_v  & J \ar [r] _ u & I . & }$$

Since $\Phi$ is a fibration there is a cartesian morphism  $\kappa
: X' \to X$ over $u$. By the cartesian property, there is a unique
morphism $\theta': Z \to X'$ over $v$ such that $\kappa
\theta'=\theta$, as in the diagram
$$\xybiglabels\xymatrix@R=2.8pc@C=2.5pc{Z \ar @{-->}[r]_{\theta'}
\ar@/^1.5pc/ [rr] ^\theta & X' \ar [r]_{\kappa} & X.}$$

Now, suppose $\phi: Y \to X$ over $u:J \to I$ satisfies $\phi
\psi= \theta$, as in the diagram:
$$\xybiglabels\xymatrix@R=2.8pc@C=2.5pc{Z \ar @{-->}[dr]_{\theta'}\ar
[r]_\psi
\ar@/^1.5pc/ [rr] ^\theta & Y \ar @{..>}[r] _ \phi & X \\
& X' \ar [ur]_{\kappa}& }$$

By the given property of $\psi$ there is a unique morphism $\psi':
Y \to X'$ in $\X_J$ such that $\psi'\psi=\theta'$. By the
cartesian property of $\kappa$, there is a unique morphism $\phi'$
in $\X_J$ such that $\kappa \phi'=\phi$. Then $$
\kappa\psi'\psi=\kappa \theta'=\theta=\phi\psi=  \kappa \phi'\psi
.$$ By the cartesian property of $\kappa$, and since
$\psi'\psi,\;\phi' \psi $ are over $uv$, we have $\psi'\psi=\phi'
\psi$. By the given property of $\psi$, and since $\phi', \;\psi'$
are in $\X_J$, we have $\phi'=\psi'$. So $\phi= \kappa\psi'$, and
this proves uniqueness.

But we have already checked that $\kappa \psi' \psi=\theta$, so we
are done.
\end{proof}

The following Proposition allows us to prove  that a fibration is
also a cofibration by constructing the adjoints $u_*$ of $u^*$ for
every $u$.

\begin{prop} \label{prop:induc-cocart}
Let $\Phi : \X \to \B$ be a fibration of categories.  Let $u: J \to
I$ have reindexing functor $u^*: \X_I \to \X_J$. Then the following
are equivalent:
\begin{enumerate}[\rm (i)] \item For all $Y \in \X_J$, there is a morphism $u_Y: Y \to
u_* Y$ which is cocartesian over $u$; \item there is a functor
$u_*: \X_J \to \X_I$ which is left adjoint to $u^*$.
\end{enumerate}
\end{prop}
\begin{proof}
That (ii) implies (i) is clear, using Proposition
\ref{prop:cocartesianweakening},  since the adjointness gives the
bijection required for the cocartesian property.

To prove that (i) implies (ii) we have to check that the
allocation $Y \mapsto u_*(Y)$ gives a functor that is adjoint to
$u^*$. As before the adjointness comes from the cocartesian
property.

We leave to the reader the check the details of the functoriality of
$u_*$.
\end{proof}

To end this section, we give a useful result on compositions.
\begin{prop}\label{prop:comp-fibrations}
The composition of fibrations, (cofibrations), is also a fibration
(cofibration).
\end{prop}
\begin{proof}
We leave this as an exercise.
\end{proof}

\section{Pushouts and cocartesian morphisms}
\label{sec:pushouts} Here is a small result which we use  in this
section and section \ref{sec:examples},  as it applies to many
examples, such as the fibration $\Ob:\Gpd \to \Set$.

\begin{prop}\label{prop:initialinfibre}
Let $\P: \X \to \B$ be a functor that has a left adjoint $\D$. Then
for each $K \in \Ob \B$, $\D(K)$ is initial in $\X_K$. In fact if
$u:K \to J$ in $\B$, then for any $X \in \X_J$ there is a unique
morphism $\eps_K : DK \to X$ over $u$.
\end{prop}
\begin{proof}
This follows immediately from the adjoint relation $\X_u(\D K,X)
\cong \B(K, \P X)$ for all $ X \in \Ob \X_J$.
\end{proof}

Special cases of cocartesian morphisms are used in
\cite{B2006,BH78:sec,BH81:col}, and we review these in section
\ref{sec:examples}. A  construction which arises naturally from the
various Higher Homotopy van Kampen theorems  is given a general
setting as follows:

\begin{thm} \label{prop:induced-pushout}
Let $\P: \X \to \B$ be a fibration of categories which has a left
adjoint $\D$. Suppose that $\X$ admits pushouts. Let $v: K \to J$ be
a morphism in $\B$, and let $Z \in \X_K$.  Then a cocartesian
lifting $\psi: Z \to Y$ of $v$ is given precisely by the pushout in
$\X$:
\begin{equation}\label{equ:pushcart}
 \xybiglabels \vcenter{\xymatrix@=3pc {\D (K) \ar [d] _{\eps_K} \ar [r] ^{\D (v)}
 & \D (J) \ar [d] ^{{\eps_J}} \\
 Z \ar[r] _{\psi}   & Y         }} \tag{*}
\end{equation}

\end{thm}
\begin{proof}
Suppose first that diagram (*)  is a pushout in $\X$. Let $u: J \to
I$ in $\B$ and let  $\theta: Z \to X$ satisfy $\P(\theta)=uv$, so
that $\P(X)=I$. Let $f: D(J) \to X$ be the adjoint of $u:J \to \P
(X)$.
\begin{equation}
 \xybiglabels \vcenter{\xymatrix {\D (K) \ar [d] _{\eps_K} \ar [r] ^{\D (v)}
 & \D (J) \ar [d] _{{\eps_J}} \ar[ddr]^f & \\
 Z  \ar [drr] _\theta \ar[r] ^{\psi}   & Y \ar@{..>}[dr] |(0.4)\phi &   \\
 && X \\
 K \ar[r]_v& J\ar[r]_u & I  }} \tag{**}
\end{equation}

Then $\P(f D(v))=uv= \P(\theta \eps_K)$ and so by Proposition
\ref{prop:initialinfibre}, $f D(v)=\theta \eps_K$. The pushout
property implies there is a unique $\phi: Y \to X$ such that $\phi
\psi = \theta$ and $\phi \eps_J= f$. This last condition gives
$\P(\phi) = u$ since $u = \P(f)=\P(\phi \eps_J)= \P(\phi) \id_J =
\P(\phi)$.

For the converse, we suppose given $f: \D(J) \to X$ and $\theta: Z
\to X$ such that $\theta {\eps_K}= f \D(v)$. Then $\P(\theta) = uv$
and so there is a cocartesian lifting ${\phi}: Y \to X $ of $u$. The
additional condition ${\phi} \eps_J= f$ is immediate by Proposition
\ref{prop:initialinfibre}.
\end{proof}

\begin{cor}
Let $\P: \X \to \B$ be a fibration which has a left adjoint and
suppose that $\X$ admits pushouts. Then $\P$ is also a cofibration.
\end{cor}
In view of the construction of hierarchical homotopical invariants
as colimits from the HHvKTs in \cite{BH81:col} and \cite{BL87}, the
following is worth recording, as a consequence of Theorem
\ref{thm:preservecolimfibres}.
\begin{thm}\label{thm:compcolim-fibred} Let $\Phi : \X \to \B$ be
fibred and cofibred. Assume $\B$ and all fibres $\X_I$ are
cocomplete. Let $T: \C \to \X$ be a functor from a small connected
category. Then $\colim T$ exists and may be calculated as
follows:\begin{enumerate}[\rm (i)]\item  First calculate $I=\colim
(\Phi   T)$, with cocone $\gamma: \Phi   T \Rightarrow I$; \item for
each $c \in \C$ choose cocartesian morphisms $\gamma'(c): T(c) \to
X(c)$, over $\gamma(c)$ where $X(c)\in \X_I$; \item make $c \mapsto
X(c)$ into a functor $X: \C \to \X_I$, so that $\gamma'$ becomes a
natural transformation $\gamma': T \Rightarrow X$;
\item form $Y=\colim X \in \X_I$ with cocone $\mu: X \Rightarrow Y$.
\end{enumerate}
Then $Y$ with $\mu \gamma': T \Rightarrow Y$ is $\colim T$.
\end{thm}
\begin{proof}
We first explain how to make $X$ into a functor.

We will in stages build up the following diagram:
\begin{equation}  \xybiglabels       \vcenter{   \xymatrix@C=3pc{& T(c) \ar @/^2pc/ [rrrr] ^\eta
\ar [dr]_{T(f)} \ar [rr]_{\gamma'(c)}
 && X(c) \ar @{-->}[d] ^{X(f)} \ar [r] ^{\mu (c)}&Y \ar [d] ^1  \ar @{..>} [r]_\tau & Z \\
 \ar[d] ^\Phi && T(c') \ar [r] _{\gamma'(c')} & X(c')\ar [r] _{\mu(c')}&Y\ar @{..>} [r]_{\tau'} &Z \\
 & K \ar [r]_{\Phi T (f)} & J \ar [r]_{\gamma(c')} & I\ar [r] _1 &I \ar [r] _w&H }}
\end{equation}
Let $f:c \to c'$ be a morphism in $\C$,  $K=\Phi T (c), J= \Phi T
(c')$.  By cocartesianness of $\gamma'(c)$, there is a unique
vertical morphism $X(f): X(c) \to X(c')$ such that $X(f)\gamma'(c)=
\gamma'(c') T(f)$. It is easy to check, again using cocartesianness,
that if further $g:c' \to c''$, then $X(gf)=X(g)X(f)$, and $X(1) =
1$. So $X$ is a functor and the above diagram shows that $\gamma'$
becomes a natural transformation $T \Rightarrow X$.

Let $\eta: T \Rightarrow Z$ be a natural transformation to a
constant functor $Z$, and let $\Phi( Z)=H$. Since $I=\colim (\Phi
T)$, there is a unique morphism $w: I \to H$ such that $w\gamma
=\Phi( \eta)$.

By the cocartesian property of $\gamma'$, there is a natural
transformation $\eta': X \Rightarrow Z$ such that $\eta' \gamma' =
\eta$.

Since $Y$ is also a colimit in $\X$ of $X$, we obtain a  morphism
$\tau: Y \to Z$ in $\X$ such that $\tau \mu = \eta'$. Then $\tau \mu
\gamma' = \eta' \gamma' = \eta$.

Let $\tau': Y \to Z$ be another morphism such that $\tau'\mu \gamma'
= \eta$. Then $\Phi(\tau)=\Phi(\tau')=w$, since $I$ is a colimit.
Again by cocartesianness, $\tau' \mu =\tau \mu$. By the colimit
property of $Y$, $\tau=\tau'$.
\end{proof}
This with Theorem \ref{thm:compcolim-fibred} shows how to compute
colimits of connected diagrams in the examples we discuss in
sections \ref{sec:groupoidsinduced} to \ref{sec:crossedsquares}, and
in all of which a van Kampen type theorem is available giving
colimits of algebraic  data for some glued topological data.

\begin{cor}
  Let $\Phi: \X \to \B$ be a functor satisfying the assumptions of
  theorem \ref{thm:compcolim-fibred}. Then $\X$ is connected
  complete, i.e. admits colimits of all connected diagrams. 
\end{cor}

\section{Groupoids bifibred over sets}\label{sec:groupoidsinduced}
We have already seen in Example \ref{example:Ob} that the  functor
$\Ob: \Gpds \to \Sets$ is a fibration. It also has a left adjoint
$\D$ assigning to  a set $I$ the discrete groupoid on $I$, and a
right adjoint assigning to a set $I$ the codiscrete groupoid on $I$.

It follows from general theorems on algebraic theories that the
category $\Gpds$ is cocomplete, and in particular admits pushouts,
and so it follows from previous results that $\Ob: \Gpds \to \Sets$
is also a cofibration. A  construction of the cocartesian liftings
of $u: I \to J$ for $G$ a groupoid over $I$ is given in terms of
words, generalising the construction of free groups and free
products of groups, in \cite{H71,B2006}.  In these references the
cocartesian lifting of $u$ to $G$ is called a {\it universal
morphism}, and is written $u_*: G \to U_u(G)$. This construction  is
of interest as it yields a normal form for the elements of $U_u(G)$,
and hence $u_*$ is injective on the set of non-identity elements of
$G$.

A homotopical application  of this cocartesian lifting is the
following theorem on the fundamental groupoid. It shows how
identification of points of a discrete subset of a space can lead to
`identifications of the objects' of the fundamental groupoid:
\begin{thm}\label{thm:vktlifting}
Let $(X,A)$ be a pair of spaces such that $A$ is discrete and the
inclusion $A \to X$ is a closed cofibration. Let $f:A \to B$ be a
function to a discrete space $B$. Then the induced morphism \[
\pi_1(X,A) \to \pi_1(B\cup _f X,B)
\]is the cocartesian lifting of $f$.
\end{thm}
This theorem immediately gives the fundamental group of the circle
$S^1$  as the infinite cyclic group $\C$, since $S^1$ is obtained
from the unit interval $[0,1]$ by identifying $0$ and $1$, as shown
in the Introduction in diagram (\ref{eq:circle1}). The theorem  is a
translation of \cite[9.2.1]{B2006}, where the words `universal
morphism' are used instead of `cocartesian lifting'. Section 8.2 of
\cite{B2006} shows how free groupoids on directed graphs are
obtained by a generalisation of this example.

The calculation of colimits in a fibre $\Gpd_I$ is similar to that
in the category of groups, since both categories are protomodular,
\cite{borcbourn}. Thus a colimit is calculated as a quotient of a
coproduct, where quotients are themselves obtained by factoring by a
normal subgroupoid. Quotients are discussed in \cite{H71,B2006}.

Theorem \ref{thm:compcolim-fibred} now shows how to compute general
colimits of groupoids.

We refer again to \cite{H71,B2006} for further developments and
applications of the algebra of groupoids; we generalise some aspects
to modules, crossed modules and crossed complexes in the next
sections.

\section{Groupoid modules bifibred over groupoids}\label{sec:modules}
Modules over groupoids are a useful generalisation of modules over
groups, and also form part of the basic structure of crossed
complexes. Homotopy groups $\pi_n(X;X_0), n \geq 2,$ of a space $X$
with a set $X_0$ of base points form a module over the  fundamental
groupoid $\pi_1(X,X_0)$, as do the homotopy groups $\pi_n(Y,X:X_0),
n \geq 3,$ of a pair $(Y,X)$.
\begin{Def} \label{defn:modulegpd}
A {\it module over a groupoid} \index{module!over a groupoid} is a
pair $(M, G),$ where $G$ is a groupoid with set of objects $G_0$,
$M$ is a totally disconnected abelian groupoid with the same set of
objects as $G$, and with a  given  action of $G$ on $M$. Thus $M$
comes with a {\it target function} $t:M \to G_0$, and each
$M(x)=t^{-1}(x), x \in G_0$, has the structure of Abelian group. The
$G$-action is given by  a family of maps
$$M(x) \times G(x,y) \to M(y)$$ for all $x,y \in G_0$. These maps are denoted by
$(m,p) \mapsto m^p$ and satisfy the usual properties, i.e. $m^1=m$,
$(m^p)^{p'}=m^{(pp')}$ and $(m+m')^p=m^p+{m'}^p$, whenever these are
defined. In particular, any $p \in G(x,y)$ induces an isomorphism $m
\mapsto m^p $ from $M(x)$ to $M(y)$.

A {\it morphism of modules} \index{morphism!module}
\index{module!morphism} is a pair $(\theta , f ):(M, G) \rightarrow
(N,H),$ where $f : G \rightarrow H$ and $\theta: M \to N$ are
morphisms of groupoids and preserve the action. {That is,} $\theta$
is given by a family of group morphisms
$$\theta(x): M(x) \rightarrow N(f(x))$$ for all $x \in G_{0}$
satisfying $\theta(y) (m^p ) = (\theta(x) (m))^{f (p)} ,$ for all $p
\in G(x,y), m \in M(x)$.

This  defines the category  \index{Mod@{$\Mod$}}$\Mod$ having
modules over groupoids as objects and the morphisms of modules as
morphisms. If $(M,G)$ is a module, then $(M,G)_0$ is defined to be
$G_0$. \end{Def}

We have a forgetful functor $\Phi_\M: \Mod \to \Gpds$ in which $(M,
G) \mapsto G$.
\begin{prop}
The forgetful functor $\Phi_\M: \Mod \to \Gpds$ has a left adjoint
and is fibred and  cofibred.
\end{prop}
\begin{proof}
The left adjoint of $\Phi_\M$  assigns to a groupoid $G$ the  module
written $0 \to G$ which has only the trivial group over each $x \in
G_0$.

Next, we give the pullback construction to prove that  $\Phi_\M$ is
fibred. This is entirely analogous to the group case, but taking
account of the geometry of the groupoid.

So let $v: G \to H$ be a morphism of groupoids, and let $(N,H)$ be a
module. We define $(M,G)=v^*(N,H)$ as follows. For $x \in G_0$ we
set $M(x)=\{x\}\times N(vx)$ with addition given by that in $N(vx)$.
The operation is given by $(x,n)^p=(y,n^{vp})$ for $ p\in G(x,y)$.

The construction of $N= v_*(M,G)$ for $(M,G)$ a $G$-module is as
follows.

For $y \in H_0$ we let $N(y)$ be the abelian group generated by
pairs $(m,q)$ with $m \in M, q \in H$, and $t(q)=y, s(q)=v(t(m))$,
so that $N(y)=0 $ if no such pairs  exist. The operation of $H$ on
$N$ is given by $(m,q)^{q'}= (m,qq')$, addition is
$(m,q)+(m',q)=(m+m',q)$ and the relations imposed are
$(m^p,q)=(m,v(p)q)$ when these make sense. The cocartesian morphism
over $v$ is given by $\psi: m \mapsto (m,1_{vt(m)})$.
\end{proof}
\begin{rem}\label{rem:attachingcells}
We will discuss the relation between a module over a groupoid and
the restriction to the vertex groups in section
\ref{sec:crossedcomplex} in the general context of crossed
complexes. However it is useful to give the  general situation of
many base points to describe the relative homotopy group
$\pi_n(X,A,a_0)$ when $X$ is obtained from $A$ by adding $n$-cells
at various base points. The natural invariant to consider is then
$\pi_n(X,A,A_0)$ where $A_0$ is an appropriate set of base points.
\end{rem}
We now describe free modules over groupoids in terms of the
inducing construction. The interest of this is two fold. Firstly,
induced modules arise in homotopy theory from a HHvKT, and we get
new proofs of results on free modules in homotopy theory.
Secondly, this indicates the power of the HHvKT since it gives new
results.
\begin{Def}
Let $Q$ be a groupoid. A {\it free basis} for a module $(N,Q)$
consists of a pair of functions $t_B:B \to Q_0$, $i: B \to N$ such
that $t_Ni=t_B$ and with the universal property that if $(L,Q)$ is a
module and $f:B \to L$ is a function such that $t_Lf=t_N$ then there
is a unique $Q$-module morphism $\phi: N \to L$ such that $\phi i
=f$.
\end{Def}

\begin{prop}\label{prop:inducedgivesfreemod}
Let $B$ be an indexing set, and $Q$ a groupoid. The  free
 $Q$-module $(FM(t),Q)$ on $t: B \to Q_0$ may be described as  the $Q$-module
induced by $t: B \to Q$ from the  discrete free module $\Z
B=(\Z\times B,B)$ on $B$, where $B$ denotes also the discrete
groupoid on $B$.
\end{prop}
\begin{proof}This is a direct verification of the universal property.
\end{proof}
\begin{rem}
Proposition \ref{prop:cocartesianweakening} shows that the universal
property for a free module can also be expressed in terms morphisms
of modules $(FM(t),Q) \to (L,R)$.
\end{rem}

\section{Crossed modules bifibred over
groupoids}\label{sec:xmodGpd} Out homotopical example here is the
family of second relative homotopy groups of a pair of spaces with
many base points.

A {\it crossed module over a groupoid}, \cite{BH81:algcub}, consists
first of a morphism of groupoids $\mu:M \to P$ of groupoids with the
same set $P_0$ of objects such that $\mu$ is the identity on
objects, and $M$ is a family of groups $M(x), x \in P_0$; second,
there is an action of $P$ on the family of groups $M$ so that if $m
\in M(x)$ and $p\in P(x,y)$ then $m^p \in M(y)$. This action must
satisfy the usual axioms for an action with the additional
properties: \\
CM1) $\mu(m^p)= p^{-1}\mu(m)p$, and \\
CM2) $m^{-1}nm= n^{\mu m}$\\
for all $p \in P$, $m,n \in M$ such that the equations make sense.
These form the objects of the category $\XMod$ in which a morphism
is a commutative square of morphisms of groups
$$\xybiglabels\xymatrix{M\ar
[d]_\mu \ar [r] ^g & N \ar [d]^\nu \\
P \ar [r]_f& Q }$$which preserve the action in the sense that
$g(m^p)=(gm)^{fp}$ whenever this makes sense.

The category $\XMod$  is equivalent to the well known category
$\2-Gpd$ of 2-groupoids, \cite{BH81:inf}.  However the advantages of
$\XMod$ over 2-groupoids are:
\begin{itemize}  \item crossed modules are closer to the classical
invariants of relative homotopy groups;
\item the notion of freeness is clearer in $\XMod$ and models a
standard topological situation, that of attaching 1- and 2-cells;
\item the category $\XMod$ has a monoidal closed structure which helps to define a notion of
homotopy; these constructions  are simpler to describe in detail
than those for 2-groupoids, and they extend to all dimensions.
\end{itemize}

We have a forgetful functor $\Phi_1: \XMod \to \Gpds$ which sends
$(M \to P) \mapsto P$. Our first main result is:
\begin{prop}
The forgetful functor $\Phi_1: \XMod \to \Gpds$ is fibred and has a
left adjoint.
\end{prop}
\begin{proof}
The left adjoint of $\Phi_1$  assigns to a groupoid $P$ the crossed
module $0 \to P$ which has only the trivial group over each $x \in
P_0$.

Next, we give the pullback construction to prove that  $\Phi_1$ is
fibred. So let $f: P \to Q$ be a morphism of groupoids, and let
$\nu: N \to Q$ be a crossed module. We define $M=\nu^*(N)$ as
follows.

For $x \in P_0$ we set $M(x)$ to be the subgroup of $ P(x) \times
N(fx)$ of elements $(p,n)$ such that $fp=\nu n$.  If $p_1 \in
P(x,x'), n \in N(fx)$ we set $(p,n)^{p_1}= (p_1^{-1}pp_1,
n^{f(p_1)})$, and let $\mu: (p,n)\mapsto p $.  We leave the reader
to verify that this gives a crossed module, and that the morphism
$(p,n) \mapsto n$ is cartesian.
\end{proof}

The following result in the case of crossed modules of groups
appeared in \cite{BH78:sec}, described in terms of the crossed
module $\partial: u_*(M) \to Q$ {\it induced} from the crossed
module $\mu: M \to P$ by a morphism $u: P \to Q$.
\begin{prop}
The forgetful functor $\Phi_1: \XMod \to \Gpds$ is cofibred.
\end{prop}
\begin{proof}
We prove  this  by a direct construction.

Let $\mu: M \to P$ be a crossed module, and let $f:P \to Q$ be a
morphism of groupoids. The construction of $N= f_*(M)$ and of
$\partial: N \to Q$ requires just care to the geometry of the
partial action  in addition to the construction for the group case
initiated in \cite{BH78:sec} and pursued in \cite{BW03} and the
papers referred to there.

Let $y \in Q_0$. If there is no $q \in Q$ from a point of $f(P_0)$
to $y$, then we set $N(y)$ to be the trivial group.

Otherwise, we define $F(y)$ to be the free group on the set of pairs
$(m,q)$ such that $m \in M(x)$ for some $x \in P_0$ and $q \in Q(f
x, y)$. If $q' \in Q(y,y')$ we set $(m,q)^{q'}= (m,qq')$. We define
$\partial': F(y) \to Q(y)$ to be $(m,q) \mapsto q^{-1}(fm)q$. This
gives a precrossed module over $\partial: F \to Q$, with function
$i: M \to F$ given by $m \mapsto (m,1)$ where if $m \in M(x)$ then
$1$ here  is the identity in $Q(fx)$.

We now wish to change the function $i: M \to F$ to make it an
operator morphism. For this,  factor $F$ out by the relations
\begin{align*}
  (m,q)(m',q)&= (mm',q),\\
  (m^p,q) &= (m,(fp)q),
\end{align*}
whenever these are defined, to give a projection $F \to F'$ and $i':
M \to F'$. As in the group case, we have to check that $\partial': F
\to Q$ induces $\partial'': F' \to H$ making this a precrossed
module. To make this a crossed module involves factoring out Peiffer
elements, whose theory is as for the group case in  \cite{BHu}. This
gives a crossed module morphism $(\phi,f): (M,P) \to (N, Q)$ which
is cocartesian.
\end{proof}
We recall the algebraic origin of free crossed modules, but in the
groupoid context.

Let $P$ be a groupoid, with source and target functions written
$s,t: P \to P_0$. A subgroupoid $N$ of $P$ is said to be {\it
normal} in $P$, written $N \lhd P$, if $N$ is wide in $P$, i.e.
$N_0 = P_0$, and for all $x,y \in P_0$ and $a \in P(x,y)$,
$a^{-1}N(x)a = N(y)$. If $N$ is also totally intransitive, i.e.
$N(x,y)= \emptyset$ when $x \ne y$, as we now assume, then the
quotient groupoid $P /N$ is easy to define. (It may also be
defined in general but we will need only this case.)

Now suppose given a family $R(x), x \in P_0$ of subsets of $P(x)$.
Then the normaliser $N_P(R )$ of $R$ is well defined as the
smallest normal subgroupoid of $P$ containing all the sets $R(x)$.
Note that the elements of $N_P(R )$ are all {\it consequences} of
$R$ in $P$, i.e. all well defined products of the form
\begin{equation} c= (r_1^{\epsilon_1})^{a_1}\ldots
(r_n^{\epsilon_n})^{a_n}, \quad \epsilon_i=\pm 1, a_i \in P, n
\geq 0
\end{equation}and where $b^a$ denotes $a^{-1}b a$. The quotient $P/N_P(R )$ is
also written $P/R$, and called the {\it quotient} of $P$ by the
relations $r=1, r \in R$.

As in group theory, we need also to allow for repeated relations. So
we suppose given a set $R$ and a function $\omega: R \to P$ such
that $s \omega = t \omega= \beta$, say. This `base point function',
saying where the relations are placed, is a useful part of the
information.

We now wish to obtain `syzygies' by replacing the normal subgroupoid
by a `free object' on the relations $\omega:R \to P$. As in the
group case, this is done using {\it free crossed modules}.
\begin{rem}There is a subtle reason for this use of crossed modules.
A normal subgroupoid $N$ of $P$ (as defined above) gives a quotient
object $P/N$ in the category $\Gpd _X$ of groupoids with object set
$X=P_0$. Alternatively, $N$ defines a congruence on $P$, which is a
particular kind of equivalence relation. Now an equivalence relation
is in general a particular kind of subobject of a product, but in
this case, we must take the product in the category $\Gpd _X$. As a
generalisation of this, one should take a groupoid object in the
category  $\Gpd _X$. Since these totally disconnected normal
subgroupoids determine equivalence relations on each $P(x,y)$ which
are congruences, it seems clear that a groupoid object internal to
$\Gpd_X$ is equivalent to a 2-groupoid with object set $X$.   \qed
\end{rem}
\begin{Def}
A {\it free basis } for a crossed module $\partial: C \to P$ over a
groupoid $P$ is a set  $R$, function $\beta:R \to P_0$ and function
$i: R \to C$ such that $i(r)  \in  C(\beta r), r \in R$, with the
universal  property that if $\mu: M \to P$ is a crossed module and
$\theta: R \to M$ a function over the identity on $P_0$ such that
$\mu \theta=
\partial i$, then there is a unique morphism of crossed $P$-modules
$\phi: C\to M$ such that $\phi i=\theta$. \qed
\end{Def} \begin{example} Let $R$ be
a set and $\beta: R \to P_0$ a function. Let $\id: F_1(R)\to F_2(R)$
be the identity crossed module on two copies of  $F(R)$, the
disjoint union of copies $\C(r)$ of the infinite cyclic group $\C$
with generator $c_r \in \C(r)$. Thus $F_2(R)$ is a totally
intransitive groupoid with object set $R$. Let $i: R \to F_1(R)$ be
the function $ r \mapsto c_r$. Let $\beta: R \to R$ be the identity
function. Then $\id: F_1(R) \to F_2(R)$ is the free crossed module
on $i$. The verification of this is simple from the diagram
$$\xybiglabels \vcenter{\xymatrix{R \ar@/^1.5pc/[rr]^\theta \ar [r]_i & F_1(R) \ar [d] _\id \ar @{-->}[r]_-f& M \ar [d]^\mu \\
& F_2(R) \ar [r] _\id& F_2(R)}}$$The morphism $f$ simply maps the
generator $c_r$ to $\theta r$. \qed
\end{example}
\begin{prop}
Let $R$ be a set, and $\mu: M \to P$ a crossed module over the
groupoid $P$.  Let $\beta: R \to P_0$ be a function. Then the
functions $i:R \to M$ such that $s\mu=t \mu =\beta $ are bijective
with the crossed module morphisms $(f,g)$
$$\xybiglabels \xymatrix{F_1(R) \ar [d] _\id \ar [r]^-f& M \ar [d]^\mu \\
F_2(R) \ar [r] _-g& P}$$ such that $sg= \beta$.

Further, the free crossed module $\partial: C(\omega) \to P$ on a
function $\omega: R \to P$ such that $s\omega=t\omega=\beta$ is
determined as the crossed module induced from $\id: F_1(R) \to
F_2(R)$ by the extension of $\omega$ to the groupoid morphism
$F_2(R) \to P$.
\end{prop}
\begin{proof}
The first part is clear since $g=\mu f$ and $f$ and $i$ are related
by $f(c_r)=i(r), r \in R$.

The second part follows from the first part and the universal
property of induced crossed modules as shown in the following
diagram: \begin{equation*} \xybiglabels \vcenter{\xymatrix{F_1(R)
\ar
[d] _\id \ar @/^1pc/ [rr]^ \theta \ar [r]_-f & C(\omega) \ar @{-->}[r]_\phi \ar [d]^\partial & M \ar [d]^\mu \\
F_2(R) \ar [r] _-g & P \ar [r] _= & P }} \qed
\end{equation*}
\end{proof}

\section{Crossed complexes and an HHvKT}
\label{sec:crossedcomplex} Crossed complexes are  analogous to chain
complexes but also generalise groupoids to all dimensions and with
their base points and operations relate dimensions 0,1 and $n$. The
structure and axioms for a crossed complex are those universally
satisfied by the main topological example, the {\it fundamental
crossed complex} $\Pi X_*$ of a filtered space $X_*$, where $(\Pi
X_*)_1$ is the fundamental groupoid $\pi_1(X_1,X_0)$ and for $n \geq
2$ $(\Pi X_*)_n$ is the family of relative homotopy groups
$\pi_n(X_n, X_{n-1},x_0)$ for all $x_0 \in X_0$, with associated
operations of the fundamental groupoid and boundaries.

Crossed complexes fit into our scheme of algebraic structures over a
range of dimensions satisfying a HHvKT in that the fundamental
crossed complex functor
$$\Pi: (\text{filtered spaces})\to (\text{crossed complexes})$$
preserves certain colimits. We state a precise version below.

A crossed complex $C$ is in part a sequence  of the form
$$\xybiglabels\xymatrix{\cdots \ar[r]& C_n \ar[r]^-{\delta_n} & C_{n-1} \ar[r]^-{\delta_{n-1}}&
\cdots \cdots \ar[r]^-{\delta_3}& C_2 \ar[r]^-{\delta_2} & C_1 }$$
where all the $C_n, n \geq 1$, are groupoids over $C_0$. Here
$\delta_2: C_2 \to C_1$ is a crossed module and for $n\geq 3$
$(C_n,C_1)$ is a module. The further axioms are that $\delta_n$ is
an operator morphism for $n \geq 2$ and that $\delta_2 c$ operates
trivially on $C_n$ for $c \in C_2$ and $n \geq 3$. We assume the
basic facts on crossed complexes as surveyed in for example
\cite{Brown-grenoble,Brown-fields}. The category of crossed
complexes is written $\Crs$. A full exposition of the theory of
crossed complexes will be given in \cite{BHS}.

To state   the {\it Higher Homotopy van Kampen Theorem for relative
homotopy groups}, namely Theorem C of \cite[Section 5]{BH81:col}, we
need the following definition:

\begin{Def}
A filtered space $X_*$ is said to be {\it connected}
\index{connected!filtered space} \index{filtered spaces!connected}
if the following conditions hold for each $n \geqslant 0:$ \\
- $\phi (X_*, 0) :$ If $r > 0,$ the map $\pi_0 X_0 \rightarrow \pi_0
X_r,$ induced by inclusion, is surjective;  i.e. $X_0$ meets
all path connected components of all stages of the filtration $X_*$. \\
 -$\phi (X_*, n)$(for $n \geqslant
1$): If $r > n$ and $x \in X_0$, then $ \pi_n (X_r , X_n , x) =
0.$\hfill $\Box$ \end{Def}

\begin{thm}[Higher Homotopy van Kampen Theorem]  \label{thm:HHvKTpushout}
Let $X_*$ be a filtered space and suppose:
\begin{enumerate}[\rm (i)]  \item  $X$ is the union of the interiors of subspaces $U$, $V$;
\item the  filtrations $U_*, V_*$ and $W_*$, formed by intersection with $X_*$,
and where $W=U \cap V$, are connected filtrations.
Then\end{enumerate}
{\rm (Conn)} the filtration  $X_*$ is connected, and \\
{\rm (Pushout)} the following   diagram of  morphisms of crossed
complexes induced by inclusions
$$\xymatrix{\Pi W_* \ar [d] \ar [r] & \Pi U_* \ar [d] \\
\ar [r]\Pi V_* & \Pi X_*} $$  is a pushout of crossed complexes.
\end{thm}
\begin{rem} The connectivity conclusion is significant, but not as important
as the algebraic conclusion. This theorem is proved without recourse
to traditional methods of algebraic topology such as homology and
simplicial approximation. Indeed, the implications in dimension 2
are  in general nonabelian and so  unreachable by the traditional
abelian methods. Instead the theorem is proved using the
construction of $\rho X_*$, the cubical homotopy $\omega$--groupoid
of the filtered space $X_*$, defined  in dimension $n$ to be the set
of filter homotopy classes of maps $I^n_* \to X_*$. The properties
of this construction enable the proof of the 1-dimensional theorem
van Kampen theorem to be generalised to higher dimensions, and the
theorem on crossed complexes is deduced using a non trivial
equivalence between the two constructions.

The paper \cite{BH81:col} also proves a more general theorem, in
which arbitrary unions lead to a coequaliser rather than a pushout.
The paper also  assumes a $J\,_0$ condition on the filtered spaces;
but this can be relaxed by the  refined definition of making filter
homotopies of maps $I^n_* \to X_*$ to be rel vertices, as has been
advertised in \cite{BH91}.\QED
\end{rem}

\begin{rem}
Colimits of crossed complexes may be computed from the colimits of
the groupoids, crossed modules and modules from which they are
constituted.
\end{rem}

\begin{rem}\label{rem:pushoutwarning}
A warning has to be given that some of the algebra is not as
straightforward as that in traditional homological and homotopical
algebra. For example in an abelian category, a pushout of the form
$$\xybiglabels \xymatrix{A \ar [r] \ar [d]_i & 0 \ar [d] \\
B \ar [r]_p & C }$$ is equivalent to an exact sequence
$$ A \labto{i} B \labto{p} C \to 0.$$

However in the category $\Mod_*$ of modules over groups a pushout of
the form
$$\xybiglabels \xymatrix{(M,G) \ar [r] \ar [d]_i & (0,1) \ar [d] \\
(N,H)  \ar [r]_p & (P,K) }$$ is equivalent to a pair of an exact
sequence of groups
$$ G \labto{i} H \labto{p} K \to 1,$$
and  an exact sequence of induced modules over $K$
$$ M \otimes_{\Z G} \Z K \labto{i} N \otimes_{\Z H} \Z K  \labto{p} P \to 0. $$
This shows that pushouts give much more information in our case, but
also shows that handling the information may be more difficult, and
that misinterpretations could lead to false conjectures or proofs.
\hfill $\Box$
\end{rem}

\begin{prop}
The truncation functor $tr^1:\Crs \to \Gpd$, $C \mapsto C_1$, is a
bifibration.
\end{prop}
\begin{proof} The previous results give the constructions on modules
and crossed modules. The functoriality of these constructions give
the construction of the boundary maps, and the axioms for all these
follow.
\end{proof}

We will also need for later applications (Proposition
\ref{prop:retractoffree}) the relation of a connected crossed
complex to the full, reduced (single vertex) crossed complex it
contains, analogous to the well known relation of  a connected
groupoid to any of its  vertex groups.

Recall that a {\it codiscrete groupoid} $T$ is one on which $T(x,y)$
is a singleton for all objects $x,y \in T_0$. This is called a {\it
tree groupoid} in \cite{B2006}. Similarly, a {\it codiscrete crossed
complex} $T$ is one in which the groupoid $T_1$ is codiscrete and
which is trivial in higher dimensions.

We follow similar conventions for crossed complexes as for groupoids
in \cite{B2006}. Thus if $D$ and $E$ are crossed complexes, and
$S=D_0 \cap E_0$ then by the {\it free product} $D*E$ we mean the
crossed complex given by the pushout of crossed complexes
\begin{equation}\label{equ:freeproduct}
  \vcenter{\xymatrix{S \ar [r] \ar [d] & E \ar [d]^j \\
  D\ar [r]_-i & D*E }}
\end{equation}
where the set $S$ is identified with the discrete crossed complex
which is $S$ in dimension 0 and trivial in higher dimensions.  The
following result is analogous to and indeed includes standard facts
for groupoids (cf. \cite[6.7.3, 8.1.5]{B2006}).

\begin{prop}\label{prop:crossedcomplex:conn}Let $C$ be a connected
crossed complex, let $x_0 \in C_0$ and let $T$ be codiscrete wide
subcrossed complex of $C$.  Let $C(x_0)$  be the subcrossed complex
of $C$ at the base point $x_0$. Then  the natural morphism
$$\phi: C(x_0) * T \to C$$
determined by the inclusions is an isomorphism, and  $T$ determines
a strong deformation retraction
$$r: C \to C(x_0). $$

Further, if $f:C \to D$ is a morphism of crossed complexes which is
the identity on $C_0 \to D_0$ then we can find a retraction $s: D
\to D(x_0)$ giving rise to a pushout square
\begin{equation} \label{eq:retractsquare}
\xybiglabels  \vcenter{\xymatrix{C \ar [d] _f \ar [r] ^-r& C(x_0) \ar [d] ^{f'}\\
  D \ar [r] _-s & D(x_0)}}
\end{equation}
in which $f'$ is the restriction of $f$.
\end{prop}
\noindent {\bf Proof} Let $i:C(x_0)\to C, \, j: T \to C$ be the
inclusions. We verify the universal property of the free product.
Let $\a: C(x_0) \to E, \b: T \to E$ be morphisms of crossed
complexes agreeing on $x_0$. Suppose $g:C \to E$ satisfies $gi=\a,
gj=\b$. Then $g$ is determined on $C_0$. Let $c \in C_1(x,y)$. Then
\begin{align*}
  c&=(\tau x)((\tau x) ^{-1}c (\tau y))(\tau
y)^{-1}\tag{*}\\[-2ex]
\intertext{and so } gc&=g(\tau x)g((\tau x) ^{-1}c (\tau y))g(\tau
y)^{-1}\\&= \b(\tau x) \a ((\tau x) ^{-1}c (\tau y))\b(\tau y)^{-1}.\\[-2ex]
\intertext{If $c \in C_n(x), n \geq 2$, then} c&=(c^{\tau x})^{(\tau
x)^{-1}}\tag{**}\\[-2ex]
\intertext{and so } g(c)&= \a(c^{\tau x})^{\b(\tau x)^{-1}}.
\end{align*}
This proves uniqueness of $g$, and conversely one checks that this
formula defines a morphism $g$ as required.

In effect, equations (*) and (**) give  for the elements of $C$
normal forms in terms of elements of $C(x_0)$ and of $T$.

This isomorphism and the constant map $T \to \{x_0\}$ determine the
strong deformation retraction $r: C \to C(x_0)$.

The retraction $s$ is defined by the elements $f\tau(x), x \in C_0$,
and then the diagram (\ref{eq:retractsquare}) is a pushout since it
is a retract of the pushout square
\begin{equation*} \label{eq:retractsquare2}
\xybiglabels  \vcenter{\xymatrix{C \ar [d] _f \ar [r] ^1& C \ar [d] ^{f}\\
  D \ar [r] _1 & D}} \tag*{$\Box$}
\end{equation*}

\section{Homotopical excision and induced constructions}
\label{sec:examples} We now interpret the HHvKT (Theorem
\ref{thm:HHvKTpushout}) when the filtrations have essentially just
two stages.

\begin{Def} \label{def:basedpair-filtration}
By a {\it based pair} $X_\bullet=(X,X_1;X_0)$ of spaces we mean a
pair $(X,X_1)$ of spaces together with a set $X_0 \subseteq X_1$
of {\it base points}. For  such a based pair and  $n \geq 2$ we
have an associated filtered space $X_\bullet^{[n]}$ which is $X_0$
in dimension $0$, $X_1$ in dimensions $1 \leq i < n$ and $X$ for
dimensions $i \geq n$. We write $\Pi_n X_\bullet$ for the crossed
complex $\Pi X_\bullet^{[n]}$. This crossed complex is trivial in
all dimensions $\neq 1,n$, and in dimension $n$ is the family of
relative homotopy groups $\pi_n(X,X_1,x_0)$ for $x_0 \in X_0$,
considered as a module (crossed module if $n=2$) over the
fundamental groupoid $\pi_1(X_1,X_0)$. Colimits of such crossed
complexes are equivalent to colimits of the corresponding module
or crossed module. \QED
\end{Def}

The following is clear.

\begin{prop} If $X_\bullet=(X,X_1;X_0)$ is a based pair, and $n \geq 2$
then the associated filtered space $X_\bullet^{[n]}$ is connected
if and only if the based pair $X_\bullet $ is {\it
$(n-1)$-connected}, i.e. if:
\begin{itemize}
\item $X_0$ meets each path component of $X_1$ and of $X$; \item each
path in $X$ joining points of $X_0$ is deformable in $X$ rel end
points  to a path in $X_1$; and \item $\pi_r(X,X_1,x_0)=0$ for all
$x_0 \in X_0$ and $1<r<n$.
\end{itemize} The last
part of the condition is of course vacuous if $n=2$. \hfill $\Box$
\end{prop}

With this in mind, Theorem \ref{thm:HHvKTpushout} may be restated
as:
\begin{thm} \label{thm:HHvKTpushoutbased}
Let $(X,X_1;X_0)$ be a based pair of  spaces and suppose:
\begin{enumerate}[\rm (i)]  \item  $X$ is the union of the interiors of subspaces $U$, $V$;
\item the  based pairs $(U, U_1; U_0), (V, V_1; V_0)$ and $(W,
W_1; W_0)$ formed by intersection with $(X,X_1;X_0)$, and where
$W=U \cap V$, are $(n-1)$-connected. Then
\end{enumerate}
{\rm (Conn)} $(X,X_1;X_0)$ is $(n-1)$-connected, and \\
{\rm (Pushout)} the following   diagram of  morphisms
$$\xymatrix{\Pi_n (W, W_1; W_0) \ar [d] \ar [r] & \Pi_n (U, U_1; U_0) \ar [d] \\
\ar [r]\Pi_n (V, V_1; V_0) & \Pi_n (X,X_1;X_0) } $$  is a pushout
of modules if $n \geq 3$ and of  crossed modules if $n=2$.
\end{thm}
\begin{proof}The important point is the equivalence between colimits of crossed
complexes which for a given $n>2$ have $C_i=0$ for $i \ne 0,1,n$ or
which are trivial in dimensions $>2$, and  colimits of the
corresponding modules or crossed modules.
\end{proof}
\begin{rem}
It is not easy to see for  Theorem \ref{thm:HHvKTpushoutbased} a
direct proof  in terms of modules or crossed modules, since one
needs the intermediate structure between 0 and $n$ to use the
connectivity conditions. The cubical homotopy $\omega$-groupoid with
connections $\rho X_*$ is found convenient for this inductive
process  in \cite{BH81:col}. There are two reasons for this: cubical
methods are convenient for constructing homotopies, and also for
{\it algebraic inverses to subdivision}. \QED
\end{rem}

We now concentrate on excision, since this gives rise to cocartesian
morphisms and so induced modules and crossed modules.

\begin{thm}[Homotopy Excision Theorem] \label{thm:homexcision-cocartesian}
Let the  topological space $X$ be the union of the interiors of sets
$U,V$, and let $W= U \cap V$. Let $n \geq 2$. Let $W_0 \subseteq U_0
\subseteq U$ be such that the based pair $(V,W;W_0)$ is
$(n-1)$-connected and $U_0$ meets each path component of $U$. Then
$(X,U;U_0)$ is $(n-1)$-connected, and the morphism of modules
(crossed if $n=2$)
$$ \pi_n(V,W;W_0) \to  \pi_n(X,U;U_0)$$
induced by inclusions is cocartesian over the morphism of
fundamental groupoids
$$ \pi_1(W,W_0) \to \pi_1(U,U_0)$$
induced by inclusion.
\end{thm}
\begin{proof}
We deduce this excision theorem from the  pushout theorem
\ref{thm:HHvKTpushoutbased}, applied to the  based pair $(X,U;U_0)$
and the following diagram of morphisms :\begin{equation}
\xymatrix{\Pi_n(W,W;W_0) \ar [r]
\ar [d] & \Pi_n(U,U;U_0)\ar [d] \\
\Pi_n(V,W;W_0) \ar [r] & \Pi_n(X,U;U_0) }
\end{equation}This is a pushout of modules if $n \geq 3$ and of crossed modules if
$n=2$, by Theorem \ref{thm:HHvKTpushoutbased}. However
$$\Pi_n(W,W;W_0)= (0,\pi_1(W,W_0)), \qquad \Pi_n(U,U;U_0)= (0,\pi_1(U,U_0)). $$
So the theorem follows from Theorem \ref{prop:induced-pushout} and
our discussion of the examples of modules and crossed modules.
\end{proof}

\begin{cor}[Homotopical excision for an adjunction] \label{cor:homexc-adjunction}
Let $i:W \to V$ be a closed cofibration and $f: W \to U$ a map.
Let $W_0$ be a subset of $W$ meeting each path component of $W$
and $V$, and let $U_0$ be a subset of $U$ meeting reach path
component of $U$ and such that $f(W_0) \subseteq U_0$. Suppose
that the based pair $(V,W)$ is $(n-1)$-connected. Let $X=U \cup _f
V$. Then the based pair $(X,U)$ is $(n-1)$-connected and the
induced  morphism of modules (crossed if $n=2$)
$$ \pi_n(V,W;W_0) \to  \pi_n(X,U;U_0)$$
is cocartesian over the induced morphism of fundamental groupoids
$$ \pi_1(W,W_0) \to \pi_1(U,U_0). $$
\end{cor}
\begin{proof}
This follows from Theorem \ref{thm:homexcision-cocartesian} using
mapping cylinders in a similar manner to the proof of a
corresponding result for the fundamental groupoid
\cite[9.1.2]{B2006}. That is, we form the mapping cylinder
$Y=M(f)\cup W$. The closed cofibration assumption ensures that the
projection from $Y$ to $X=U\cup _f V$ is a homotopy equivalence.
\end{proof}
\begin{rem}
These methods were used in \cite{BH78:sec}. \QED
\end{rem}
\begin{cor}[Attaching $n$-cells]\label{cor:freeonattaching}
Let the space $Y$ be obtained from the space $X$ by attaching
$n$-cells, $n \geq 2$,  at a set of base points $A$ of $X$, so that
$Y=X \cup_{f_\lambda}e^n_\lambda, \lambda \in \Lambda$, where
$f_\lambda:(S^{n-1},0) \to (X,A)$. Then $\pi_n(Y,X;A)$ is isomorphic
to the free $\pi_1(X,A)$-module (crossed if $n=2$) on the
characteristic maps of the $n$-cells.
\end{cor}

\begin{rem}The previous corollary for $n=2$ was a theorem of J.H.C.
Whitehead. An account of Whitehead's proof is given in \cite{B80}.
There are several other proofs in the literature but none give the
more general homotopical excision result, theorem
\ref{thm:homexcision-cocartesian}. \QED
\end{rem}
\begin{example}\label{ex:intro} We can now explain the example in
the Introduction. That $S^n$ is $(n-1)$-connected and
$\pi_n(S^n,0)\cong \Z$ follows by induction in the usual way from
the homotopical excision theorem and the calculation
$\pi_1(S^1,0)\cong \Z$ by the groupoid van Kampen theorem. Applying
the HET to writing $S^n \vee [0,1]$ as a union of $S^n$ and $[0,1]$
we get that $\pi_n(S^n\vee [0,1], [0,1];\{0,1\})$ is the free
$\pi_1([0,1],\{0,1\})$-module on one generator. Again applying the
HET but now identifying $0,1$ we get that $\pi_1(S^n \vee S^1,0)$ is
the free $\pi_1(S^1,0)$-module in one generator. \QED
\end{example}
\begin{cor}[Relative Hurewicz Theorem] \label{cor:RHT}Let $A \to X$ be a closed
cofibration and suppose $A$ is path connected and $(X,A)$ is
$(n-1)$-connected. Then $X\cup CA$ is $(n-1)$-connected and
$\pi_n(X\cup CA,x)$ is isomorphic to $\pi_n(X,A,x)$ factored by the
action of $\pi_1(A,x)$.
\end{cor}
We now point out that a generalisation of a famous result of Hopf,
\cite{hopf-secondbettigroup}, is a  corollary of the relative
Hurewicz theorem.  The following for $n=2$ is part of Hopf's result.
The algebraic description of $H_2(G)$ which he gives for $G$ a group
is shown in \cite{BH78:sec} to follow from the HHvKT.
\begin{prop}[Hopf's theorem]
Let $(V,A)$ be a pair of pointed spaces such that:
\begin{enumerate}[\rm (i)] \item $\pi_i(A)=0$ for $1 < i <n$;
\item $\pi_i(V)=0$ for $1 < i \leqslant n$; \item the inclusion $A
\to V$ induces an isomorphism on fundamental groups.
\end{enumerate}
Then the pair $(V,A)$ is $n$-connected, and the inclusion $A \to V$
induces an epimorphism $H_nA \to H_n V$ whose kernel consists of
spherical elements, i.e. of the  image of $\pi_n A$ under the
Hurewicz morphism $\omega_n: \pi_n(A) \to H_n(A)$.\end{prop}
\begin{proof} That $(V,A)$ is $n$-connected follows immediately  from the homotopy
exact sequence of the pair $(V,A)$ up to $\pi_n(V)$. We now consider
the next part of the exact homotopy sequence and its relation to the
homology exact sequence as shown in the commutative diagram:
\begin{equation*}
\xybiglabels \xymatrix{\pi_{n+1}(V,A) \ar [d] _{\omega_{n+1}} \ar
[r] ^-{\partial} & \pi_n(A)
\ar [d]_{\omega_n} \ar [r] & \pi_n(V) \ar [d] \ar [r] & \pi_n(V,A) \ar [d] \\
H_{n+1}(V,A)  \ar [r]_-{\partial'}  & H_n(A)  \ar [r]_{i_*} & H_n(V)
\ar [r] & H_n(V,A) }
\end{equation*}
The Relative Hurewicz Theorem implies that $H_n(V,A)=0$, and that
$\omega_{n+1}$ is surjective. Also   $\partial$ in the top row is
surjective, since $\pi_n(V)=0$. It follows easily that the sequence
$\pi_n(A) \to H_n(A) \to H_n(V) \to 0$ is exact.
\end{proof}

There is a nice generalisation of Corollary
\ref{cor:freeonattaching}, which so far  has been proved only as a
deduction from a HHvKT.
\begin{cor}[Attaching cones]\label{cor:attachcone} Let $A$ be a space and let
$S$ be a set consisting of one point in each path component of $A$.
By $CA$, the cone on $A$,   we mean the union of cones on each path
component of $A$. Let $f: A \to X$ be a map, and let $S'$ be the
image of $S$ by $f$. Then $\pi_2(X\cup CA,X;S')$ is isomorphic to
the $\pi_1(X,S')$-crossed module induced from the identity crossed
module $\pi_1(A,S) \to \pi_1(A,S)$ by the induced morphism $f_*:
\pi_1(A,S) \to \pi_1(X,S')$.
\end{cor}
The paper \cite{BW03} uses this result to give explicit calculations
for the crossed modules representing the homotopy 2-types of certain
mapping cones.

We now explain the relevance to  free crossed modules of Proposition
\ref{prop:crossedcomplex:conn}, leaving the module and other cases
to the reader.
\begin{prop}\label{prop:retractoffree}Let $X$ be a path connected
space with base point $a$, and let $Y= X \cup
_{f_\lambda}\{e^2_\lambda\}$ be obtained by attaching cells by means
of pointed maps $f_\lambda:(S^1,0) \to (X,a_\lambda)$, determining
elements $x_\lambda\in \pi_1(X,a_\lambda)$, $\lambda \in \Lambda$.
Let $A=\{a\}\cup \{a_\lambda\mid \lambda \in \Lambda\}$. Let $T$ be
a tree groupoid in $\pi_1(X,A)$ determining a retraction $r:
\pi_1(X,A)\to \pi_1(X,a)$. Then $\pi_2(Y,X,a)$ is isomorphic to the
free crossed $\pi_1(X,a)$-module on the elements $r(x_\lambda)$.
\end{prop}
\begin{proof} We consider the following diagram:
\begin{equation}
  \xybiglabels \vcenter{\xymatrix{(0\to \bigsqcup_\lambda \Z ) \ar [d] \ar [r] &
  (0 \to P) \ar [d] \ar [r]^r &
  (0 \to P(a)) \ar [d] \\
  (\bigsqcup_\lambda \Z \to \bigsqcup_\lambda \Z ) \ar [r] & (C(\Lambda) \to P) \ar [r]_s &(F \to P(a))  }}
\end{equation}
The left hand square is the pushout defining the free crossed module
$C(\Lambda)\to P$ as an induced crossed module. The right hand
square is the  special case of crossed modules of the retraction of
Proposition \ref{prop:crossedcomplex:conn}, and so is also a
pushout. Hence the composite square is a pushout. Hence the crossed
module $F \to P(a)$ is the free crossed module as described.
\end{proof}

\begin{rem}
An examination of Whitehead's  paper \cite{W41}, and the exposition
of part of it in \cite{B80}, shows that the geometrical side of the
last proposition is intrinsic to his approach. Of course a good
proportion of Whitehead's work was devoted to extending ideas of
combinatorial group theory to higher dimensions in combinatorial
homotopy theory. The argument here is that this extension requires
{\it combinatorial groupoid theory} for good modelling of the
geometry. \QED
\end{rem}


\section{Crossed squares and triad homotopy groups}\label{sec:crossedsquares}
In  this section we give a brief sketch  of the theory of triad
homotopy groups, including the exact sequence  relating them to
homotopical excision, and show that the  third triad group forms
part of a {\it crossed square} which, as an algebraic structure
with links over several dimensions, in this case dimensions 1,2,3,
fits
 our criteria for a HHvKT. Finally we indicate a bifibration from
crossed squares, so leading to the notion of {\it induced crossed
square}, which is relevant to a triadic Hurewicz theorem in
dimension 3.

A {\it triad of spaces} $(X:A,B;x)$ consists of a pointed space
$(X,x)$ and two pointed subspaces $(A,x),(B,x)$. Then
$\pi_n(X:A,B;x)$ is defined for $n \geq 2$ as the set of homotopy
classes of maps  $$(I^n:
\partial^-_1I^n,\partial^-_2I^n;J^{n-1}_{1,2}) \to (X:A,B;x)$$where
$J^{n-1}_{1,2}$ denotes the union of the faces of $I^n$ other than
$\partial^-_1I^n,\partial^-_2I^n$. For $n \geq 3$ this set obtains
a group structure, using the direction 3, say, which is Abelian
for $n \geq 4$. Further there is an exact sequence
\begin{equation} \to \pi_{n+1}(X:A,B;x) \to \pi_n(A,C,x)
\labto{\epsilon} \pi_n(X,B,x) \to\pi_{n}(X:A,B;x) \to
\end{equation}where $C=A \cap B$, and $\epsilon$ is the excision
map.  It was the fact that these groups measure the failure of
excision that was their main interest. However they do not shed
light on the above Homotopical Excision Theorem
\ref{thm:homexcision-cocartesian}.

The third triad homotopy group fits into a diagram of possibly
non-Abelian groups \begin{equation}\label{eq:triadcrossed square}
  \Pi(X;A,B,x):= \vcenter{\xymatrix{\pi_3(X;A,B,x) \ar [d]\ar [r] &\pi_2(B,C,x) \ar [d] \\
  \pi_2(A,C,x) \ar [r] & \pi_1(C,x)  }}
\end{equation}
in which $\pi_1(C,x)$ operates on the other groups and there is also
a function $$\pi_2(A,C,x) \times \pi_2(B,C,x) \to \pi_3(X:A,B;x)$$
known as the generalised Whitehead product.

This diagram has structure and properties which   are known as those
of a {\it crossed square}, \cite{LGW,BL87}, explained below, and so
this gives a homotopical functor
\begin{equation} \Pi: (\text{based triads}) \to (\text{crossed
squares}).
\end{equation}

A { \it crossed square} is a commutative diagram of morphisms of
groups
\begin{equation}\label{eq:crossedsquare}
  \vcenter{\xymatrix{L \ar [d] _{\lambda '} \ar [r] ^\lambda & M \ar [d] ^\mu \\
  N \ar [r] _\nu & P }}
\end{equation}
together with  left actions of $P$ on $L,M,N$ and a function $h: M
\times N \to L$ satisfying a number of axioms which we do not give
in full here. Suffice it to say that the morphisms in the square
preserve the action of $P$, which acts on itself by conjugation;
$M,N$ act on each other and on $L$ via $P$; $\lambda, \lambda',
\mu, \nu$ and $\mu \lambda$ are crossed modules; and $h$ satisfies
axioms reminiscent of commutator rules, summarised by saying it is
a  {\it biderivation}. Morphisms of crossed squares are defined in
the obvious way, giving a category $\XSq$ of crossed squares.

Let $\XMod^2$ be the category of pairs of crossed modules $\mu: M
\to P, \nu:N\to P$ (with $P$  and $\mu, \nu$ variable), and with the
obvious notion of morphism. There is a forgetful functor $\Phi: \XSq
\to\XMod^2$. This functor has a right adjoint $\D$ which completes
the pair $\mu : M \to P, \nu : N \to P$ with $L= M \times_PN$ and
$\lambda, \lambda'$ given by the projections and $h: M \times N \to
L$ given by $h(m,n)= (^nmm\io, n\, ^mn\io),m \in M, n \in N $. More
interestingly, it has a left adjoint which to the above pair of
crossed $P$-modules yields the `universal crossed square'
\begin{equation}
\vcenter{\xymatrix{M \otimes N \ar [d] \ar [r] & N \ar[d] ^\nu\\
M \ar [r] _ \mu & P}}
\end{equation}
where $M \otimes N$, as defined in \cite{BL87},  is the nonabelian
tensor product of groups which act on each other.

Then $\Phi$ is a fibration of categories and also a cofibration.
Thus we have a notion of {\it induced crossed square}, which
according to Proposition \ref{prop:induced-pushout} is given by a
pushout  of the form
$$\vcenter{\xymatrix@C=5pc@R=4pc{\ar [d]_{\begin{pmatrix} u &
1 \\1 & 1 \end{pmatrix}}{\begin{pmatrix} M \otimes N & N \\M & P
\end{pmatrix} } \ar [r]^-{\begin{pmatrix} \alpha \otimes \beta &
\beta\\ \alpha & \gamma
\end{pmatrix} } & {\begin{pmatrix} R \otimes S & S \\ R & Q
\end{pmatrix}} \ar [d]^{\begin{pmatrix} v &
1 \\1 & 1 \end{pmatrix}} \\ {\begin{pmatrix} L & N \\M & P
\end{pmatrix} } \ar [r]_-{\begin{pmatrix} \delta &
\beta\\ \alpha & \gamma
\end{pmatrix} } & {\begin{pmatrix} T & S \\R & Q \end{pmatrix}} }}
$$
in the category of crossed squares, given morphisms $(\alpha,
\gamma): (M \to P) \to (R \to Q), (\beta,\gamma): (N \to P) \to (S
\to Q)$ of crossed modules.

The functor $\Pi$ is exploited in \cite{BL87} for an HHvKT
implying some calculations of the non-Abelian group $\pi_3(X\colon
A,B;x)$\footnote{Earlier results had used homological methods to
obtain some Abelian values.}. The applications  are developed  in
\cite{BL2} for a triadic Hurewicz Theorem, and for the notion of
free crossed square, both based on `induced crossed squares'. Free
crossed squares are exploited in \cite{Ellis-3-type} for homotopy
type calculations.

In fact the HHvKT  works in all dimensions and in the more general
setting of $n$-cubes of spaces, although not in a many base point
situation. For a recent application, see \cite{ellis-mikhailov}.

\end{document}